\documentclass{article}

\newcommand{\smalltitre}{Petit titre}

\usepackage{amsmath}

\def\aut{\Large \centerline{\auts}}

\def\titre{\sl}

\def\auts{\sl}

\newcommand{\tit}{

\begin{center}

\begin{minipage}{12cm}

\begin{center}

\bf\Large\titre

\end{center}

\end{minipage}

\end{center}

\rm\normalsize

\aut

\addcontentsline{toc}{part}{\textsc{\titre}\\\rm\auts}\normalsize\vspace{1cm}}

\def \add#1#2{\noindent\bf{#1}\\\rm\noindent#2\par}

\renewenvironment{abstract}{\bf ABSTRACT: \sl}{\rm}

\newtheorem{lem}{Lemma}[section]

\newtheorem{pro}[lem]{Proposition}

\newtheorem{theo}[lem]{Theorem}

\usepackage{amssymb}

\begin{document}

\renewcommand{\titre}{Probabilistic Analysis for Randomized Game Tree Evaluation} 
\renewcommand{\smalltitre}{Randomized Game Tree Evaluation}
\renewcommand{\auts}{T\"amur Ali Khan and Ralph Neininger}
\tit

\begin{abstract}
We give a probabilistic analysis for the randomized game tree evaluation algorithm of
Snir. We first show that there exists an input such that the running time, measured as the
number of external nodes read by the algorithm,  on that input is maximal in stochastic
order among all
possible inputs. For this worst case input we
identify the exact expectation of the number of external nodes read by the algorithm,
give the asymptotic order of the variance including the leading constant,
provide a
limit law for an appropriate normalization as well as a tail bound
estimating large deviations. Our tail bound improves upon the exponent
of an earlier bound due to Karp and Zhang, where subgaussian tails were
shown based on an approach using multitype branching processes and
Azuma's inequality. Our approach rests on a direct, inductive estimate
of the moment generating function.
\end{abstract}

\index{game tree}
\index{randomized algorithm}
\index{multitype branching process}
\index{minimax tree}
\index{limit law}
\index{large deviations}
\index{tail bound}
\index{recursive algorithm}

\section{Introduction}
In this note we analyze the performance of the randomized
algorithm to evaluate Boolean decision trees proposed by Snir (1985). Given is
a complete binary tree of height $2k$, $k\ge 1$, where the root (at depth $0$) is
labeled $\wedge$ as are all internal nodes with even depth, all internal
nodes
with odd depth 
are labeled $\vee$.
The $n=2^{2k}$ external nodes are labeled either $0$ or $1$ and the objective is to
calculate the value of the root. For each node its value is given as the value of the
operation labeled at that node applied to the values of its children. The cost
 for
evaluating the Boolean decision tree is measured as the number of external nodes
read by the algorithm.

Snir proposed and analyzed the following randomized algorithm to evaluate
a Boolean decision tree: At each node one chooses randomly (with probability $1/2$)
one of its children and calculates its value recursively. If the result allows to identify
the value of the node (that is a $0$ for a $\wedge$-labeled node and a $1$ for
 a $\vee$-labeled node, respectively) one is done, otherwise also the other child's
value has to be calculated recursively in order to obtain the value of the node.
Applying this to the root of the tree yields the value of the Boolean decision tree.

The advantage of this algorithm over any deterministic algorithm is that for any input
at the external nodes its expected cost is sublinear in $n$, whereas any deterministic
algorithm has linear worst case cost. More precisely, Saks and Wigderson (1986)
obtained that the maximum expected cost is of the order $\Theta(n^\alpha)$ with
$\alpha=\log_2((1+\sqrt{33})/4)\doteq 0.753$ and showed that this is also a
lower bound on the maximum expected cost for any other randomized algorithm
to evaluate a Boolean decision tree; see also Motwani and Raghavan (1995, Chapter 2)
for an account on this subject. Further analysis was given by Karp and
Zhang (1995). For certain regular inputs at the external nodes the cost
of the algorithm can be represented via  2-type Galton-Watson
processes. Karp and Zhang showed that the normalized cost has
subgaussian tails. That argument was based on Azuma's inequality.

We denote the input of $0$'s and $1$'s at the external nodes as a vector $v\in\{0,1\}^n$
and the number of external nodes read by the algorithm on input $v$ by $C(v)$. We
will see subsequently that for particular $v^\star\in\{0,1\}^n$ not only the expectation
of the cost of the algorithm is maximized, i.e., 
$\mathbb{E}\, C(v^\star)=\max_{v\in\{0,1\}^n}\mathbb{E}\, C(v)$,
but also that $C(v^\star)$ is maximal in stochastic order, $C(v)\preceq C(v^\star)$ for all
$v\in\{0,1\}^n$. Here, $X\preceq Y$ for random variables $X,Y$ denotes that the
corresponding distribution functions $F_X, F_Y$ satisfy $F_X(x)\ge F_Y(x)$ for all
$x\in \mathbb{R}$, or, equivalently, that there are realizations $X',Y'$
of the distributions ${\cal L}(X), {\cal L}(Y)$ of $X,Y$
on a joint probability space such that we pointwise have $X'\le Y'$.

From this perspective it is reasonable to consider $C(v^\star)$ as
 the worst case complexity of the randomized algorithm and to analyze its asymptotic
probabilistic behavior.
Our  results for the exact mean of $C(v^\star)$, the asymptotic growth of its
variance including the evaluation of the leading constant, a limit law
for $C(v^\star)$ after normalization as $k\to\infty$ together with an
explicit tail estimate are
 based on a recursive description of the problem. Since $v^\star$ is
a  regular input in the sense of Karp and Zhang, also their 2-type
Galton-Watson approach applies. 

Our main finding is an improvement
of the tail bound $\exp(-\mathrm{const}\; t^2)$ for $t>0$, to
$\exp(-\mathrm{const}\; t^\kappa)$, with
$1<\kappa<1/(1-\alpha)\doteq 4.06$, see Theorem 3.6.
This is based on a direct, inductive estimate of the moment
generating function. Our approach is also applicable to any
regular input as well as to other related problems.

The paper is organized as follows: In section 2 we explain, how a
worst case input $v^\star$ is obtained. Section 3 contains the
statements
of the results. In sections 4 and 5 the $2$-type branching process of
Karp and Zhang (1995) is recalled and the recursive description of the
quantities, that our analysis is based on, is introduced. Section 6
contains the proofs of our results and section 7 has extensions to 
$m$-ary Boolean decision trees.

\section{Worst case input}
In this section we explain how a worst case input $v^\star$ is constructed.
We first have a look at the case $k=1$ and $v\in\{0,1\}^4$ such that the decision tree
is evaluated  to $1$ at the root. Clearly both children of the root have
to lead to an evaluation of $1$.
Now each pair of external nodes attached to the children needs to have at least one value
$1$. Note that the algorithm reads in both pairs of external nodes until it finds the first
one. Hence there will in total be read two $1$'s no matter how
$v\in\{0,1\}^4$ is drawn among the choices
that lead to an evaluation of $1$ for the decision tree. Clearly, to
maximize the number of $0$'s being read
we choose in each pair of external nodes one $0$ and one $1$. Then both $0$'s are being
read independently with probability $1/2$. Hence,  $v_1=(0,1,0,1)$ stochastically maximizes
 $C(v)$ for all   $v\in\{0,1\}^4$ such
 that the decision tree evaluates $1$, see Figure 1.

Analogously look at the case $k=1$ and $v\in\{0,1\}^4$ such that the decision tree
is evaluated  to $0$. Clearly, one child of the root has to have the value $0$, whose
external nodes attached need to have both values $0$. If we choose also value $0$
for the other child of the root, we are lead to $v=(0,0,0,0)$, and the algorithm reads
exactly 2 external nodes with values both $0$. Therefore, to stochastically maximize
$C(v)$ we choose the second child of the root with value $1$ and again its external nodes
attached with values $0$ and $1$. Then, $v_0=(0,0,0,1)$ stochastically maximizes
 $C(v)$ for all   $v\in\{0,1\}^4$ for which the decision tree evaluates to $0$, see Figure 1.

Since we have $C(v_0)\preceq C(v_1)$, it follows that $v^\star=(0,1,0,1)$ is a choice with
$C(v)\preceq C(v^\star)$ for all  $v\in\{0,1\}^4$. For general $k\ge 2$ a corresponding
$v^\star=v^\star(k)$ can recursively be constructed from  $v^\star(k-1)$
as follows: 
Each component
$0$ in    $v^\star(k-1)$ is replaced by the block $0,0,0,1$, whereas each $1$ is
replaced by the block $0,1,0,1$. For example, for $k=3$, this yields
\begin{eqnarray*}
v^\star&=&(0,0,0,1,0,0,0,1,0,0,0,1,0,1,0,1,\\
&&\phantom{(}0,0,0,1,0,1,0,1,0,0,0,1,0,1,0,1,\\
&&\phantom{(}0,0,0,1,0,0,0,1,0,0,0,1,0,1,0,1,\\
&&\phantom{(}0,0,0,1,0,1,0,1,0,0,0,1,0,1,0,1).
\end{eqnarray*} 
In Proposition 3.1 we show that this construction yields a $v^\star$
with $C(v)\preceq C(v^\star)$ for all  $v\in\{0,1\}^n$ and $k\ge
1$.

\centerline{\hspace{0cm}\hrulefill \hspace{0cm}}
\begin{center}
\setlength{\unitlength}{1800sp}
\begingroup\makeatletter\ifx\SetFigFont\undefined
\gdef\SetFigFont#1#2#3#4#5{
  \reset@font\fontsize{#1}{#2pt}
  \fontfamily{#3}\fontseries{#4}\fontshape{#5}
  \selectfont}
\fi\endgroup
\begin{picture}(4294,3495)(268,-3019)
{\thicklines
\put(2415, 48){\circle{812}}
\put(1230,-1225){\circle{812}}
\put(3600,-1225){\circle{812}}
\put(2012,-188){\line(-1,-1){605}}
\put(2818,-188){\line( 1,-1){605}}
\put(638,-2290){\line(2,5){288}}
\put(1815,-2290){\line( -2,5){288}}
\put(3014,-2290){\line(2,5){288}}
\put(4191,-2290){\line( -2,5){288}}
\put(1890,-80){\makebox(0,0)[lb]{\smash{\SetFigFont{14}{34.8}
{\rmdefault}{\mddefault}{\itdefault}{$\wedge$}}}}
\put(705,-1393){\makebox(0,0)[lb]{\smash{\SetFigFont{14}{34.8}
{\rmdefault}{\mddefault}{\itdefault}{$\vee$}}}}
\put(3075,-1393){\makebox(0,0)[lb]{\smash{\SetFigFont{14}{34.8}
{\rmdefault}{\mddefault}{\itdefault}{$\vee$}}}}
\put(301,-2986){\framebox(675,675){}}
\put(2677,-2986){\framebox(675,675){}}
\put(1478,-2986){\framebox(675,675){}}
\put(3854,-2986){\framebox(675,675){}}}
\put(1950,-852){\makebox(0,0)[lb]{\smash{\SetFigFont{14}{34.8}
{\rmdefault}{\mddefault}{\itdefault}{$1$}}}}
\put(775,-2125){\makebox(0,0)[lb]{\smash{\SetFigFont{14}{34.8}
{\rmdefault}{\mddefault}{\itdefault}{$1$}}}}
\put(3145,-2125){\makebox(0,0)[lb]{\smash{\SetFigFont{14}{34.8}
{\rmdefault}{\mddefault}{\itdefault}{$1$}}}}
\put(2539,-2825){\makebox(0,0)[lb]{\smash{\SetFigFont{14}{34.8}
{\rmdefault}{\mddefault}{\itdefault}{$0$}}}}
\put(1340,-2825){\makebox(0,0)[lb]{\smash{\SetFigFont{14}{34.8}
{\rmdefault}{\mddefault}{\itdefault}{$1$}}}}
\put(3716,-2825){\makebox(0,0)[lb]{\smash{\SetFigFont{14}{34.8}
{\rmdefault}{\mddefault}{\itdefault}{$1$}}}}
\put(163,-2825){\makebox(0,0)[lb]{\smash{\SetFigFont{14}{34.8}
{\rmdefault}{\mddefault}{\itdefault}{$0$}}}}
\end{picture}
\quad\quad\quad\quad
\begin{picture}(4294,3495)(268,-3019)
{\thicklines
\put(2415, 48){\circle{812}}
\put(1230,-1225){\circle{812}}
\put(3600,-1225){\circle{812}}
\put(2012,-188){\line(-1,-1){605}}
\put(2818,-188){\line( 1,-1){605}}
\put(638,-2290){\line(2,5){288}}
\put(1815,-2290){\line( -2,5){288}}
\put(3014,-2290){\line(2,5){288}}
\put(4191,-2290){\line( -2,5){288}}
\put(1890,-80){\makebox(0,0)[lb]{\smash{\SetFigFont{14}{34.8}
{\rmdefault}{\mddefault}{\itdefault}{$\wedge$}}}}
\put(705,-1393){\makebox(0,0)[lb]{\smash{\SetFigFont{14}{34.8}
{\rmdefault}{\mddefault}{\itdefault}{$\vee$}}}}
\put(3075,-1393){\makebox(0,0)[lb]{\smash{\SetFigFont{14}{34.8}
{\rmdefault}{\mddefault}{\itdefault}{$\vee$}}}}
\put(301,-2986){\framebox(675,675){}}
\put(2677,-2986){\framebox(675,675){}}
\put(1478,-2986){\framebox(675,675){}}
\put(3854,-2986){\framebox(675,675){}}}
\put(1950,-852){\makebox(0,0)[lb]{\smash{\SetFigFont{14}{34.8}
{\rmdefault}{\mddefault}{\itdefault}{$0$}}}}
\put(775,-2125){\makebox(0,0)[lb]{\smash{\SetFigFont{14}{34.8}
{\rmdefault}{\mddefault}{\itdefault}{$0$}}}}
\put(3145,-2125){\makebox(0,0)[lb]{\smash{\SetFigFont{14}{34.8}
{\rmdefault}{\mddefault}{\itdefault}{$1$}}}}
\put(2539,-2825){\makebox(0,0)[lb]{\smash{\SetFigFont{14}{34.8}
{\rmdefault}{\mddefault}{\itdefault}{$0$}}}}
\put(1340,-2825){\makebox(0,0)[lb]{\smash{\SetFigFont{14}{34.8}
{\rmdefault}{\mddefault}{\itdefault}{$0$}}}}
\put(3716,-2825){\makebox(0,0)[lb]{\smash{\SetFigFont{14}{34.8}
{\rmdefault}{\mddefault}{\itdefault}{$1$}}}}
\put(163,-2825){\makebox(0,0)[lb]{\smash{\SetFigFont{14}{34.8}
{\rmdefault}{\mddefault}{\itdefault}{$0$}}}}
\end{picture}
\end{center}
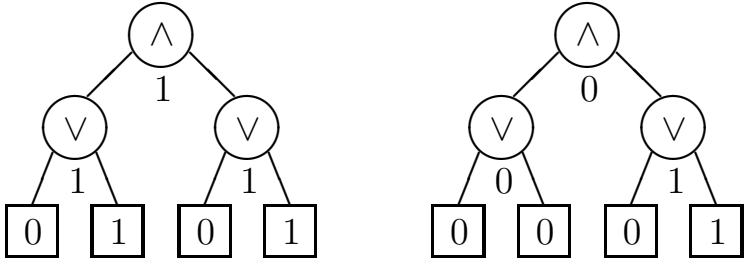
Figure 1: {\em Shown are decision trees for $k= 1$ evaluating at the root
    to $1$ and $0$, respectively,
together with a choice for the external nodes that stochastically maximizes
the number of external nodes read by the algorithm.}\\
\vspace{-3mm}
\centerline{\hspace{0cm}\hrulefill \hspace{0cm}}\vspace{3mm}\\

If we would only want to stochastically maximize the cost over all 
$v\in R_0(n)\subset \{0,1\}^n$ that evaluate
to a $0$ at the root, the same recursive construction of replacing
digits by 
corresponding blocks,
starting with $v_0=(0,0,0,1)$, yields a $v_\star \in
R_0(n)$ such that $C(v) \preceq C(v_\star)$ for all $v\in R_0(n)$.

\section{Results}
We assume that we have $n=2^{2k}$ with $k\ge 1$ and denote by
$v^\star\in\{0,1\}^n$ an input as constructed in section 2.
\begin{pro} For $v^\star\in \{0,1\}^n$ as defined in section 2 we have 
$C(v)\preceq C(v^\star)$
for all $v\in \{0,1\}^n$.
\end{pro}
The stochastic worst case behavior $C(v^\star)$ of the randomized game tree evaluation
algorithm has the following asymptotic properties: The subsequent
theorems describe the behavior of mean, variance, limit distribution,
and large deviations of $C(v^\star)$. For the mean we have:
\begin{theo}
The expectation of $C(v^\star)$ is given by 
$\mathbb{E}\, C(v^\star) = c_1 n^{\alpha} -c_2 n^{\beta}$,
with
\begin{eqnarray*}
\alpha=\log_2\frac{1+\sqrt{33}}{4}, \quad
\beta=\log_2\frac{1-\sqrt{33}}{4}   ,\quad
c_1=\frac{1}{2}+\frac{7}{2\sqrt{33}},  \quad c_2=
c_1-1.
\end{eqnarray*}
\end{theo}
We denote for sequences $(a_k), (b_k)$ by $a_k\sim b_k$ asymptotic equivalence, i.e.,
$a_k/b_k\to 1$ as $k\to \infty$. Then we have for the variance of $C(v^\star)$:
\begin{theo}
The variance of  $C(v^\star)$ satisfies asymptotically 
$\mathrm{Var}\; C(v^\star)\sim d\, n^{2\alpha}$ as $k\to\infty$,
where $d\doteq 0.0938$. The constant $d$ can also be given in closed form.
\end{theo}
For random variables $X,Y$ we denote by $X\stackrel{d}{=} Y$ equality in distribution,
i.e., ${\cal L}(X)={\cal L}(Y)$. Then we have the following limit law
for $C(v^\star)$:
\begin{theo}
For $C(v^\star)$ we have after normalization convergence in
distribution, \begin{eqnarray*} \frac{C(v^\star)}{n^{\alpha}} \longrightarrow
C, \quad k\to\infty, \end{eqnarray*} where the distribution of
$C$ is given as ${\cal L}(C)={\cal L}(G_1)$ and ${\cal L}(G)={\cal
L}(G_0,G_1)$ is characterized by $\mathbb{E}\, \|G\|^2<\infty$, $\mathbb{E}\, G=
(c_0,c_1)$, with $c_0=1/2+5/(2\sqrt{33})$, and \begin{eqnarray*} G
\stackrel{d}{=} \frac{1}{4^\alpha}\left\{ G^{(1)}+G^{(2)}+
\left[\!\begin{array}{cc} B_1B_2&0\\ 1-B_2&0
\end{array}\!\right]G^{(3)}+
 \left[\!\begin{array}{cc} 0&B_1\\ B_1&0   \end{array}\!\right]G^{(4)}\right\},
\end{eqnarray*} with $G^{(1)},\ldots,G^{(4)}$, $B_1,B_2$ independent
with ${\cal L}(G^{(r)})={\cal L}(G)$, for $r=1,\ldots,4$, and  ${\cal
L}(B_1)={\cal L}(B_2)=B(1/2)$. Here, $B(1/2)$ denotes the
Bernoulli$(1/2)$ distribution.
\end{theo}
For the estimate of large deviations we rely on Chernoff's bounding
technique. We need to follow a bivariate setting for the vector
$(C(v^\star),C(v_\star))$ as introduced in section 5. The
following bound on the moment generating function is obtained:
\begin{pro}
It exists a sequence $(Y_k)_{k\ge 0}=(Y_{k,0},Y_{k,1})_{k\ge 0}$ of bivariate random variables
with marginal distributions  ${\cal L}((C(v^\star)-\mathbb{E}\,
C(v^\star))/n^\alpha )$, ${\cal L}((C(v_\star)-\mathbb{E}\,
C(v_\star))/n^\alpha)$ such that for all $q> 1/\alpha\doteq 1.33$ there is a $K>0$ with \begin{eqnarray}\label{akne-mgf} \mathbb{E}\, \exp\langle
s, Y_k\rangle \le \exp(K\|s\|^q) \end{eqnarray} for all $s\in \mathbb{R}^2$ and
$k\ge 0$. An explicit value for $K=K_q$ is given in (\ref{akne-expl}).
\end{pro}
The bound on the moment generating function in the previous proposition 
implies a large deviation estimate via Chernoff bounds:
\begin{theo}
For all $1<\kappa< 1/(1-\alpha) \doteq 4.06$ there
exists an $L >0$ such that for any $t>0$ and $n=2^{2k}$
\begin{eqnarray}\label{akne-ldb} \mathbb{P}\left(\frac{C(v^\star)-\mathbb{E}\,
C(v^\star)}{n^\alpha}> t\right)\le \exp(-L t^\kappa). \end{eqnarray}
 An explicit value for $L$ is given
in (\ref{akne-lqt}). The same bound applies to the left tail.
\end{theo}
The approach of Karp and Zhang (1995) based on Azuma's inequality gives
the tail bound $\exp(-L't^2)$ for an explicitly known $L'$. For
$\kappa=2$ the prefactor $L=L_2$ in Theorem 3.6 can also be evaluated
and satisfies $L_2>11L'$.

\section{Karp and Zhang's $\mathbf{2}$-type branching process}
For the analysis of $C(v^\star)$ note that whenever the algorithm
has to evaluate the value of a node at a certain depth that yields
a $1$, according to the discussion of section 2, the algorithm
has to evaluate the values of two nodes of depths two levels below
that each yield a $1$, and $B_3+B_4$ nodes  of depths two levels
below that each yield a $0$, cf.~Figure 1. Here, $B_3, B_4$ are
independent Bernoulli $B(1/2)$ distributed random variables.
Analogously, when the algorithm has to evaluate the value of a
node at a certain depth that yields $0$, two levels below it has
to evaluate $B_1$ nodes yielding a $1$ and $2+B_1B_2$ nodes
yielding a $0$, where $B_1,B_2$ are independent $B(1/2)$
distributed random variables. Here, the event $\{B_1=1\}$
corresponds to the algorithm first checking the right child of the
node to be evaluated and  $\{B_2=1\}$ to first checking the left
child of that child, cf.~Figure 1. Since at each node the child
being evaluated first is independently drawn from all other
choices, this gives rise to the following 2-type Galton-Watson
branching process.

We have individuals of type $0$ and $1$ where the population of
the $k$-th generation corresponds to the number of nodes at depth
$2k$ that are read by the algorithm. We consider processes
starting either with an individual of type $1$ or type $0$ and
assume that the algorithm is applied to the worst case inputs
$v^\star$ and $v_\star$, respectively. Then we have the following
offspring distributions: An individual of type $1$ has an
offspring of $2$ individuals of type $1$ and $B_3+B_4$ individuals
of type $0$. An individual of type $0$ has an offspring of $B_1$
individuals of type $1$ and $2+B_1B_2$ individuals of type $0$. We
denote the number of individuals of type $0$ and $1$ in generation
$k$ by $(V^{(i)}_n,W^{(i)}_n)$, when starting with an individual
of type $i=0,1$, where $n=2^{2k}$. Note that for $v^\star, v_\star
\in \{0,1\}^n$ we have the representations 
\begin{eqnarray*}
 C(v^\star)\stackrel{d}{=} V^{(1)}_n+W^{(1)}_n,\quad  
C(v_\star)\stackrel{d}{=} V^{(0)}_n+W^{(0)}_n.
\end{eqnarray*} 
This is the approach of Karp and Zhang (1995) for regular
inputs like $v^\star, v_\star$. Hence, part of the analysis of
$C(v^\star)$ can be reduced to the application of the theory of
multitype branching processes; see for general reference Harris
(1963) and Athreya and Ney (1972), and for a survey on the
application of branching processes to tree structures and tree
algorithms see Devroye (1998).

However, we will also use a recursive description of the problem.
This will be given in the next section and enables to use as well results from
the probabilistic analysis of recursive algorithms by the contraction method.

\section{The recursive point of view}
It is convenient to work as well with a recursive description of
the distributions ${\cal L}(C(v_\star))$ and ${\cal
L}(C(v^\star))$. For this, we define the distributions of a
bivariate random sequence $(Z_n)=(Z_{n,0},Z_{n,1})$ for all
$n=2^{2k}$, $k\ge 0$ by $Z_1=(1,1)$ and, for $k\ge 1$,
 \begin{eqnarray*} 
Z_n \stackrel{d}{=}  Z^{(1)}_{n/4} + Z^{(2)}_{n/4} +
\left[\!\begin{array}{cc}
    B_1B_2&0\\ 1-B_2&0   \end{array}\!\right]Z^{(3)}_{n/4}
+  \left[\!\begin{array}{cc} 0&B_1\\ B_1&0
  \end{array}\!\right]Z^{(4)}_{n/4},
\end{eqnarray*} 
where $Z^{(1)}_{n/4},\ldots,Z^{(4)}_{n/4}$, $B_1,B_2$
are independent, $B_1,B_2$ are Bernoulli $B(1/2)$ distributed
and ${\cal L}(Z^{(1)}_{n/4})=\cdots={\cal L}(Z^{(4)}_{n/4})={\cal
L}(Z_{n/4})$. It can directly be checked by induction that the
marginals of $Z_n$ satisfy ${\cal L}(Z_{n,0})={\cal
L}(C(v_\star))$ and
 ${\cal L}(Z_{n,1})={\cal L}(C(v^\star))$. Note that $Z_{n,0}$ and
 $Z_{n,1}$ become dependent, firstly, since we have coupled the
 offspring
distributions using for the second component of $Z_n$ again $B_1$ and $1-B_2$ 
instead of $B_3$ and $B_4$, cf.~section 4,  and, secondly, since the first component
of $Z^{(3)}_{n/4}$ contributes to both components of $Z_n$. 
Sequences satisfying recursive
 equations as $(Z_n)$ are being dealt with in a probabilistic framework,
 the so called contraction method; see R\"osler (1991, 1992), Rachev and
R\"uschendorf (1995), R\"osler and R\"uschendorf (2001), and Neininger and
R\"uschendorf (2004).

\section{Proofs}
In this section we sketch the proofs of the results stated in section 3.\\

{\noindent {\bf Proof of Proposition 3.1:}} {\em (Sketch)\ }
We denote by $R_0(n), R_1(n)\subset \{0,1\}^n$ the sets of vectors at the
external nodes at depth $2k$ that yield an evaluation at the root of the decision tree
of value $0$ and $1$, respectively. From the discussion in section 2 we have
\begin{eqnarray*}
C(v)\preceq C(v_\star), \quad v\in R_0(n), \quad \mbox{and}\quad
C(v)\preceq C(v^\star), \quad v\in R_1(n).
\end{eqnarray*}
Hence, it remains to show that $C(v_\star)\preceq C(v^\star)$. This is
shown by 
induction on
$k\ge 1$. For $k=1$ this can directly be checked. For the step $k-1\to
k$ 
assume that we
have $C(v_\star(k-1))\preceq C(v^\star(k-1))$. It suffices to find
realizations of the quantities  $(V^{(1)}_n,W^{(1)}_n)$ and
$(V^{(0)}_n,W^{(0)}_n)$ 
on a
joint probability space with $V^{(0)}_n+W^{(0)}_n\le
V^{(1)}_n+W^{(1)}_n$ almost 
surely, $n=2^{2k}$.

For this we
use  $B,B'$,
$(V^{(i),j}_{n/4},W^{(i),j}_{n/4})$ for $i=1,2$, $j=1,\ldots,4$ being
independent for each $i=0,1$ and with  $B,B'$  Bernoulli $B(1/2)$
distributed,  ${\cal L}(V_{n/4}^{(i),j})= {\cal L}(V_{n/4}^{(i)})$, ${\cal
  L}(W_{n/4}^{(i),j})= {\cal L}(W_{n/4}^{(i)})$ for $i=1,2$ and $j=1,\ldots,4$.   By the induction
hypothesis we may assume that we have versions of these random variates with
$V^{(0),j}_{n/4}+W^{(0),j}_{n/4}
\le V^{(1),j}_{n/4}+W^{(1),j}_{n/4}$ for $j=1,\ldots, 4$. With this coupling we define
$(V^{(1)}_n,W^{(1)}_n)$ and $(V^{(0)}_n,W^{(0)}_n)$ according to the values of $B,B'$:
On  $\{B=1, B'=0\}$ we set
\begin{eqnarray*}
\left(\!\!\!\begin{array}{c} V^{(0)}_n\\ W^{(0)}_n    \end{array}\!\!\!\right) &:=&
\left(\!\!\!\begin{array}{c} V^{(0),2}_{n/4}\\ W^{(0),2}_{n/4}    \end{array}\!\!\!\right) +
\left(\!\!\!\begin{array}{c} {V}^{(0),3}_{n/4}\\ {W}^{(0),3}_{n/4}   \end{array}\!\!\!\right)+
BB' \left(\begin{array}{c}  V^{(0),4}_{n/4}\\ W^{(0),4}_{n/4}    \end{array}\!\!\!\right)
+B\left(\!\!\!\begin{array}{c} {V}^{(1),1}_{n/4}\\ {W}^{(1),1}_{n/4}
  \end{array}\!\!\!\right),
\\
\left(\!\!\!\begin{array}{c} V^{(1)}_n\\ W^{(1)}_n    \end{array}\!\!\!\right) &:=&
B\left(\!\!\!\begin{array}{c} V^{(0),3}_{n/4}\\ W^{(0),3}_{n/4}    \end{array}\!\!\!\right) +
\left(\!\!\!\begin{array}{c} {V}^{(1),1}_{n/4}\\ {W}^{(1),1}_{n/4}    \end{array}\!\!\!\right)+
B'\left(\!\!\!\begin{array}{c} V^{(0),4}_{n/4}\\ W^{(0),4}_{n/4}    \end{array}\!\!\!\right) +
\left(\!\!\!\begin{array}{c} {V}^{(1),2}_{n/4}\\ {W}^{(1),2}_{n/4}    \end{array}\!\!\!\right)
\end{eqnarray*}
and obtain $V^{(0)}_n+ W^{(0)}_n\le   V^{(1)}_n+ W^{(1)}_n$.  On the remaining sets
$\{B=0, B'=0\}$, $\{B=0, B'=1\}$, and $\{B=1, B'=1\}$ similar couplings
of  $(V^{(0)}_{n},W^{(0)}_{n})$, $(V^{(1)}_{n},W^{(1)}_{n})$ can be
defined with $V^{(0)}_n+ W^{(0)}_n\le   V^{(1)}_n+ W^{(1)}_n$.
{\;{$\rule{2mm}{2mm}$}\vskip 16pt}

{\noindent {\bf Proof of Theorem 3.2:}} {\em (Sketch)\ } Assume
that a generation has $(w_0,w_1)$ individuals of type $0$ and $1$.
Then, by the definition on the offspring distribution in section
4, the expected number of individuals in the subsequent generation
is given by 
\begin{eqnarray*} 
M\left(\!\!\!\begin{array}{c} w_0\\ w_1
\end{array}\!\!\!\right), \quad M:= \left[\!\begin{array}{cc}
9/4&1\\ 1/2&2   \end{array}\!\right]. 
\end{eqnarray*} 
Since
$C(v^\star)=C(v^\star(k))$ is the sum of the individuals at
generation $k$ for the process started with an individual of type
$1$ we obtain \begin{eqnarray*} \mathbb{E}\, C(v^\star) =
(1,1)M^k\left(\!\!\!\begin{array}{c} 0\\ 1
\end{array}\!\!\!\right). \end{eqnarray*} The matrix $M$ has the eigenvalues
$\lambda_1=(17+\sqrt{33})/8$ and $\lambda_2=(17-\sqrt{33})/8$ and
its $k$-th power can be evaluated to \begin{eqnarray*}
M^k=\frac{1}{2\sqrt{33}}\left[\!\begin{array}{cc}
(\sqrt{33}+1)\lambda_1^k+(\sqrt{33}-1)\lambda_2^k&
8(\lambda_1^k-\lambda_2^k)\\ 4(\lambda_1^k-\lambda_2^k)&
(\sqrt{33}-1)\lambda_1^k+(\sqrt{33}+1)\lambda_2^k
\end{array}\!\right]. \end{eqnarray*} From this, $\mathbb{E}\, C(v^\star)$ and various
constants needed subsequently can be read off. Note, that
$\lambda_1^k=n^\alpha$ with $\alpha$ given in Theorem 3.2 and
$n=2^{2k}$. {\;{$\rule{2mm}{2mm}$}\vskip 16pt}
Before proving Theorem 3.3 it is convenient to first prove Theorem 3.4.\\

{\noindent {\bf Proof of Theorem 3.4:}} {\em (Sketch)\ } The
2-type branching process defined in section 4 is supercritical,
nonsingular, and positive regular. Hence, a theorem of Harris
(1963) implies that \begin{eqnarray*}
\frac{1}{n^\alpha}\left(\!\!\!\begin{array}{c} V^{(1)}_n\\
W^{(1)}_n    \end{array}\!\!\!\right) \longrightarrow
Y\left(\!\!\!\begin{array}{c} \nu_1\\ \nu_2
\end{array}\!\!\!\right)  \end{eqnarray*} 
almost surely, as
$k\to\infty$, where $Y$ is a nonnegative random variable and
$(\nu_1,\nu_2)$ a deterministic vector that could also be further
specified. Thus we obtain \begin{eqnarray*}
\frac{C(v^\star)}{n^\alpha}\longrightarrow C \end{eqnarray*} in
distribution, as $k\to\infty$, with ${\cal L}(C)={\cal
L}((\nu_1+\nu_2)Y)$.

On the other hand the recursive formulation of section 5 leads
after the normalization $X_n:=Z_n/n^\alpha$ to \begin{eqnarray*} X_n \stackrel{d}{=}
\sum_{r=1}^4 A_r X^{(r)}_{n/4}, \end{eqnarray*} for $k\ge 1$, where 
$A_1=A_2=(1/4^\alpha)I_2$, with the $2\times 2$ identity matrix
$I_2$, and 
\begin{eqnarray}\label{akne-matr}
A_3=\frac{1}{4^\alpha}\left[\!\begin{array}{cc} B_1B_2&0\\ 1-B_2&0
  \end{array}\!\right],\quad
A_4=\frac{1}{4^\alpha} \left[\!\begin{array}{cc} 0&B_1\\ B_1&0
\end{array}\!\right], \end{eqnarray} 
where $X^{(1)}_{n/4},\ldots,X^{(4)}_{n/4}$, $B_1,B_2$ are 
independent with ${\cal L}(X^{(r)}_{n/4})={\cal L}(X_{n/4})$ for
$r=1,\ldots,4$ and
${\cal L}(B_1)={\cal L}(B_2)=B(1/2)$. It follows from the
contraction method that $X_n$ converges weakly and with all mixed
second moments to some $G$, that can be characterized as in
Theorem 3.4. For details, how to apply the contraction method, see
Theorem 4.1 in Neininger (2001). Thus, we have
$C(v^\star)/n^\alpha\to G_1$ in
distribution.
{\;{$\rule{2mm}{2mm}$}\vskip 16pt}

{\noindent {\bf Proof of Theorem 3.3:}} {\em (Sketch)\ } As shown
in the proof of Theorem 3.4 we have the convergence
$X_n=Z_n/n^\alpha\to G$  for all mixed second moments. This, in
particular, implies $\mathrm{Var}\,X_{n,1}\to\mathrm{Var}\, G_1$. The
variances of $G_1$  can be obtained from the distributional identity
for $G$ stated in Theorem 3.4. Then we obtain 
$\mathrm{Var}\,C(v^\star)=\mathrm{Var}(n^\alpha
X_{n,1})\sim d n^{2\alpha}$ with $d=\mathrm{Var}\, G_1$.
{\;{$\rule{2mm}{2mm}$}\vskip 16pt}

{\noindent {\bf Proof of Proposition 3.5:}}
For $Y_n=(1/n^\alpha)(Z_n-\mathbb{E}\, Z_n)$ we have marginals ${\cal
  L}(Y_{n,1})={\cal L}((C(v^\star)-\mathbb{E}\, C(v^\star))/n^\alpha)$ and  ${\cal
  L}(Y_{n,0})={\cal L}((C(v_\star)-\mathbb{E}\,
C(v_\star))/n^\alpha)$. The distributional recurrence for $Z_n$ from
section 5 implies the relation \begin{eqnarray*} Y_n \stackrel{d}{=} \sum_{r=1}^4
A_r Y^{(r)}_{n/4} +b_n,\quad k\ge 1, \end{eqnarray*} with
$Y^{(1)}_{n/4},\ldots,Y^{(4)}_{n/4}, B_1,B_2$ independent,
${\cal L}(Y^{(r)}_{n/4})={\cal L}(Y_{n/4})$,for $r=1,\ldots,4$, ${\cal
L}(B_1)={\cal L}(B_2)=B(1/2)$  and
$b_n=(1/n^\alpha)(4^\alpha\sum_{r=1}^4 (A_r\mathbb{E}\,
Z_{n/4})-\mathbb{E}\, Z_n)$. The matrices $A_r$ are given in 
(\ref{akne-matr}).

We prove the assertion by induction on $k$. For $k=0$ we have
$Y_0=0$, thus the assertion is true. Assume the assertion is true
for some $n/4=2^{2(k-1)}$. Then, conditioning on
$(A_1,\ldots,A_4,b_n)$, denoting the distribution of this vector
by $\sigma_n$, and using the induction hypothesis, we obtain \begin{eqnarray*}
\mathbb{E}\, \exp\langle s,Y_n\rangle &=&\int \exp\langle s,\beta_n\rangle
\prod_{r=1}^4 \mathbb{E}\, \exp\langle s, a_r Y_{n/4}\rangle
d\sigma_n(a_1,\ldots,a_4,\beta_n)\\
&\le& \int \exp\langle s,\beta_n\rangle \prod_{r=1}^4  \exp(K
\|a_r^Ts\|^q)
d\sigma_n(a_1,\ldots,a_4,\beta_n)\\
&\le& \int \exp\left(\langle s,\beta_n\rangle +K \|s\|^q
\sum_{r=1}^4 \|a_r\|_\mathrm{op}^q  \right)
d\sigma_n(a_1,\ldots,a_4,\beta_n)\\
&=&\mathbb{E}\, \exp(\langle s,b_n\rangle+K\|s\|^qU)\exp(K\|s\|^q), \end{eqnarray*}
with $U:=\sum_{r=1}^4
\left(\|A_r\|_\mathrm{op}^q\right)-1=4^{-\alpha
q}(2+B_1B_2+(1-B_2)+B_1)-1$ and
$\|A\|_\mathrm{op}=\sup_{\|x\|=1}\|Ax\|$ for matrices $A$. Hence, 
the proof is completed by showing \begin{eqnarray*} \sup_{k\ge1} \mathbb{E}\,
\exp(\langle s,b_n\rangle+K\|s\|^qU)\le 1, \end{eqnarray*} for some
appropriate $K>0$. We denote $\xi:=-\, \mathrm{ess\,sup}\, U  = 1-4^{1-\alpha q}$, thus
$q>1/\alpha$ implies $\xi>0$.

{\it Small $\|s\|$}: First we consider small $\|s\|$ with $\|s\|\le
c/\sup_{k\ge 1} \|b_n\|_{2,\infty}$ for some $c>0$, where
$\|b_n\|_{2,\infty}:= \|\, \| b_n\|\,\|_\infty$, the inner norm being the
Euclidean norm. Note that throughout we have $n=n(k)=2^{2k}$.
For these small $\|s\|$  we
have
\begin{eqnarray*}
 \mathbb{E}\, \exp((\langle s,b_n\rangle+K\|s\|^qU)\le \exp(-K\|s\|^q\xi)\mathbb{E}\,\exp
 \langle s,b_n\rangle
\end{eqnarray*}
and, with $\mathbb{E}\,  \langle s,b_n\rangle =0$,
\begin{eqnarray*}
\mathbb{E}\, \exp \langle s,b_n\rangle &=& \mathbb{E}\, \left[ 1+ \langle
  s,b_n\rangle+\sum_{k=2}^\infty \frac{ \langle
    s,b_n\rangle^k}{k!}\right]\\
&=& 1+ \mathbb{E}\,  \langle s,b_n\rangle^2 \sum_{k=2}^\infty \frac{ \langle
    s,b_n\rangle^{k-2}}{k!}\\
&\le& 1+\|s\|^2\mathbb{E}\, \|b_n\|^2 \sum_{k=2}^\infty \frac{c^{k-2}}{k!}\\
&=&1+ \|s\|^2\mathbb{E}\, \|b_n\|^2 \frac{e^c-1-c}{c^2}.
\end{eqnarray*}
Using $\exp(-K\|s\|^q\xi)\le 1/(1+K\|s\|^q\xi)$ and with
$\Psi(c)=(e^c-1-c)/c^2$ we obtain
\begin{eqnarray*}
 \mathbb{E}\, \exp(\langle s,b_n\rangle+K\|s\|^qU)\le \frac{1+ \|s\|^2\mathbb{E}
   \|b_n\|^2\Psi(c)}{1+K\|s\|^q\xi}.
\end{eqnarray*}
Hence, we have to choose $K$ with
\begin{eqnarray*}
K\ge \frac{\|s\|^{2-q}\Psi(c)}{\xi}\sup_{k\ge 1}\mathbb{E}\, \|b_n\|^2.
\end{eqnarray*}
With $\|s\|\le c/\sup_{k\ge 1} \|b_n\|_{2,\infty}$ a possible choice is
\begin{eqnarray*}
K= \frac{\sup_{k\ge 1}\mathbb{E}\, \|b_n\|^2}{\sup_{k\ge 1}
  \|b_n\|_{2,\infty}^{2-q}}\; \frac{\Psi_q(c)}{\xi},
\end{eqnarray*}
with $\Psi_q(c)=(e^c-1-c)/c^q$.

{\it Large $\|s\|$}: For general $s\in\mathbb{R}^2$ we have \begin{eqnarray*}
\langle s,b_n\rangle+K\|s\|^qU\le \|s\|\|b_n\|-\|s\|^qK\xi\le
\|s\| \|b_n\|_{2,\infty}-\|s\|^qK\xi, \end{eqnarray*} and this is less than
zero if \begin{eqnarray*} \|s\|^{q-1}\ge \frac{\sup_{k\ge 1}
\|b_n\|_{2,\infty}}{K\xi}= \frac{\sup_{k\ge 1}
\|b_n\|_{2,\infty}^{3-q}}{\sup_{k\ge 1}\mathbb{E}\, \|b_n\|^2\Psi_q(c)}.
\end{eqnarray*} If $\|s\|$ satisfies the latter inequality we call it large.
Thus, for large $\|s\|$ we have $\sup_{k\ge1} \mathbb{E}\, \exp(\langle
s,b_n\rangle+K\|s\|^qU)\le 1$.

In order to overlap the regions for small and large $\|s\|$ we
need \begin{eqnarray*} \Psi_1(c)\ge \frac{\sup_{k\ge 1}
  \|b_n\|_{2,\infty}^2}{\sup_{k\ge 1}\mathbb{E}\, \|b_n\|^2}.
\end{eqnarray*} The right hand side of the latter display can be evaluated
explicitly for our problem and equals
$104/77$.
Thus, this
inequality is true for, e.g., $c=1.53$. Hence, with the explicit value
\begin{eqnarray}\label{akne-expl} K:=K_q=\frac{\sup_{k\ge 1}\mathbb{E}\, \|b_n\|^2}{\sup_{k\ge
1}
  \|b_n\|_{2,\infty}^{2-q}}\;\frac{e^{1.53}-2.53}{1.53^q(1-4^{1-q\alpha})}
\end{eqnarray} the proof is completed.{\;{$\rule{2mm}{2mm}$}\vskip 16pt}

{\noindent {\bf Proof of Theorem 3.6:}} By Chernoff's bounding technique
we have, for $u>0$ and with Proposition 3.5,
\begin{eqnarray*}
\mathbb{P}\left(\frac{C(v^\star)-\mathbb{E}\, C(v^\star)}{n^\alpha}>
  t\right)&=&\mathbb{P}(\exp(uY_{n,1})>\exp(ut)))\\
&\le& \mathbb{E}\, \exp(uY_{n,1}-ut)\\
&=&\mathbb{E}\, \exp(\langle(0,u),Y_n\rangle -ut)\\
&\le& \exp(K_q u^q-ut),
\end{eqnarray*}
for all $q$, $K_q$ as in Proposition 3.5 and (\ref{akne-expl}). Minimizing over $u>0$ we obtain the
bound
\begin{eqnarray*}
\mathbb{P}\left(\frac{C(v^\star)-\mathbb{E}\, C(v^\star)}{n^\alpha}>
  t\right)\le \exp(-L t^\kappa),
\end{eqnarray*} for $1<\kappa<1/(1-\alpha)$, with \begin{eqnarray}\label{akne-lqt}
L=L_\kappa =
K_{\kappa/(\kappa-1)}^{1-\kappa}\frac{(\kappa-1)^{\kappa-1}}{\kappa^\kappa}
\end{eqnarray} and $K_{\kappa/(\kappa-1)}$ given in (\ref{akne-expl}). This
completes the tail bound.
{\;{$\rule{2mm}{2mm}$}\vskip 16pt}

\section{$\mathbf{m}$-ary Boolean decision trees}
The analysis can be carried over to the case of $m$-ary Boolean
decision trees. The algorithm visits randomly chosen children and
evaluates recursively their values until the value of the root can
be identified, the remaining children are discarded afterwards. A
worst case input $v^\star\in\{0,1\}^n$
 with $n=m^{2k}$ can be constructed similarly. Then we have similar
 results for $C(v^\star)$:
\begin{theo}
For the worst case complexity  $C(v^\star)$ of evaluating an $m$-ary
Boolean decision tree we have the following asymptotics:
\begin{eqnarray*}
\mathbb{E}\, C(v^\star) &=&c_1^{(m)} n^{\alpha_m} + c_2^{(m)} n^{\beta_m},\\
\mathrm{Var}\; C(v^\star)&\sim& d_m n^{2\alpha_m},\\
\frac{C(v^\star)}{ n^{\alpha_m}}&\to& C_m,\\
\mathbb{P}\left(\frac{C(v^\star)-\mathbb{E}\, 
C(v^\star)}{n^{\alpha_m}}>t\right)&\le& \exp(-L^{(m)}
t^\kappa), \quad t>0,
\end{eqnarray*}
with constants $c_1^{(m)}$, $\alpha_m$, $\beta_m$, $d_m$,
$L^{(m)}>0$,  $c_2^{(m)}\in \mathbb{R}$,
 and $1< \kappa< \kappa_m=1/(1-\alpha_m)$.
\end{theo}
Numerical values for $\alpha_m, d_m$ and $\kappa_m$ are 
listed in Table 1. 
The distribution of $C_m$ is
given as ${\cal L}(C_m)={\cal L}(G_1)$ and ${\cal L}(G)={\cal
  L}(G_0,G_1)$ is 
characterized
by
$\mathbb{E}\, \|G\|^2<\infty$, $\mathbb{E}\, G= (c_0^{(m)},c_1^{(m)})$ and
\begin{eqnarray*}
G &\stackrel{d}{=}& \frac{1}{m^{2\alpha_m}}\Bigg\{
\sum_{r=1}^m G^{(r)}
+\sum_{r=1}^{m-1} \left[\!\begin{array}{cc} 0&  {\bf 1}_{r}(U_0)\\   
{\bf 1}_{r}(U_0)&0\end{array}\!\right]\bar{G}^{(r)}\\
&&\phantom{ \frac{1}{m^{2\alpha_m}}\Bigg\{}
~+\sum_{r,\ell=1}^{m-1}  \left[\!\begin{array}{cc}
    {\bf 1}_{r}(U_0) {\bf 1}_{\ell}(U_r)&0\\1-{\bf 1}_{\ell}(U_r)&0
  \end{array}\!\right]G^{(r,\ell)}
        \Bigg\},
\end{eqnarray*}
with ${\cal L}(G^{(r)})={\cal L}(\bar{G}^{(r)})={\cal
  L}(G^{(r,\ell)})={\cal L}(G)$ and
$G^{(r)}$, $\bar{G}^{(r)}$, $G^{(r,\ell)}$, $U_r$ independent   with
 ${\cal L}(U_{r})=\mathrm{unif}\{0,\ldots,m-1\}$ for
all $r,\ell$. Here, we denote $ {\bf 1}_{i}(Y):= {\bf 1}_{\{i\le Y\}}$
for integer $i$ and a random variable $Y$, and we have
\begin{eqnarray*}
c_0^{(m)}=\frac{1}{2}+\frac{m+3}{2\sqrt{16m+(m-1)^2}},\quad
c_1^{(m)}=\frac{1}{2}+\frac{3m+1}{2\sqrt{16m+(m-1)^2}}.
\end{eqnarray*}
\\
\centerline{\hspace{0cm}\hrulefill \hspace{0cm}}
\begin{center}
\begin{tabular}{|l|c|c|c|c|c|c|c|}\hline
        $m$ & 2 & 3 & 4 & 5 & 6 & 7 & 8   \\ \hline
        $\alpha_m$ & 0.754 & 0.759 & 0.765 & 0.769 & 0.774 & 0.778 & 0.781   \\ \hline
        $d_m$ & 0.0938 & 0.0847 & 0.0782 & 0.0731 & 0.0689 &  0.0652 & 0.0619 \\ \hline
        $\kappa_m$ & 4.060 & 4.154 & 4.247 & 4.336 & 4.419 & 4.497 &
        4.571   \\ \hline \hline \hline \hline
        $m$ & 9 & 10 & 11 & 12 & 13 & 14 & 15   \\ \hline
        $\alpha_m$ & 0.785 & 0.788 & 0.790 & 0.793 & 0.795 & 0.798 & 0.800      \\ \hline
        $d_m$ & 0.0590 & 0.0564&  0.0541 & 0.0519 & 0.0499 & 0.0481 & 0.0464   \\ \hline
        $\kappa_m$  & 4.641 & 4.707 & 4.769 & 4.829 & 4.886 & 4.940 &
        4.993  \\ \hline \hline \hline \hline
        $m$ & 16 & 17 & 20 & 30&40&50&100 \\ \hline
        $\alpha_m$ & 0.802 & 0.804 &  0.809 &0.821&0.830&0.837& 0.856 \\ \hline
        $d_m$ & 0.0448 & 0.0433 &  0.0394 &0.0304&0.0247&0.0209& 0.0117 \\ \hline
        $\kappa_m$ & 5.043 & 5.091 &  5.226 & 5.596 & 5.885 & 6.123 & 6.928  \\ \hline
\end{tabular}\\
\end{center}
Table 1: {\em Numerical values of the  quantities $\alpha_m$, $d_m$ and
  $\kappa_m$ appearing in Theorem 7.1
  for various values of $m$.}\\
\vspace{-3mm}
\centerline{\hspace{0cm}\hrulefill \hspace{0cm}}\\

\bibliography{resume}

\add{T\"amur Ali Khan and Ralph Neininger}
{Department of Mathematics\\
J.W.~Goethe University\\
Robert-Mayer-Str.~10\\
60325 Frankfurt a.M.\\
Germany\\
\{alikhan, neiningr\}@ismi.math.uni-frankfurt.de}

\end{document}